\documentclass[12pt]{amsart}
  
\input xy
\xyoption{all}
  
\usepackage{amssymb}
\usepackage{amsmath}
\usepackage{amscd}
\usepackage{latexsym}
\usepackage{xcolor}
\usepackage{graphicx}
\usepackage{fullpage}
\usepackage{hyperref}
\begin{document}
 
\newtheorem{lemma}{Lemma}[section]
\newtheorem{prop}[lemma]{Proposition}
\newtheorem{cor}[lemma]{Corollary}
  
\newtheorem{thm}[lemma]{Theorem}
\newtheorem{Mthm}[lemma]{Main Theorem}
\newtheorem{con}{Conjecture}
\newtheorem{claim}{Claim}
\newtheorem{ques}{Question}

\theoremstyle{definition}
  
\newtheorem{rem}[lemma]{Remark}
\newtheorem{rems}[lemma]{Remarks}
\newtheorem{defi}[lemma]{Definition}
\newtheorem{ex}[lemma]{Example}

\newcommand{\C}{\mathbb C}
\newcommand{\R}{\mathbb R}
\newcommand{\Q}{\mathbb Q}
\newcommand{\Z}{\mathbb Z}
\newcommand{\N}{\mathbb N}

\title{Horospherical random graphs and lockdown strategies}

\author[I. Chatterji]{Indira Chatterji}
\address{Indira Chatterji,  LJAD, Universit\'e C\^ote d'Azur, France}
\email{indira.chatterji@math.cnrs.fr}

\author[A. Lawson]{Austin Lawson}
\address{Austin Lawson}
\email{doctorazlawson@gmail.com}

\maketitle
\tableofcontents
\begin{abstract}
Expanders are sparse graph that are strongly connected, where {\it connectivity} is quantified using eigenvalues of the adjacency matrix, and {\it sparsity} in terms of vertex degree. We give a model of random graphs and study their connectivity and sparsity. This model is a particular case of soft geometric random graphs, and allows to construct sparse graphs with good expansion properties, as well as highly clustered ones. On those graphs, we study the speed at which random walks spread in the graph, and visit all vertices. As an illustration, we build a model for mainland France and study the spread of random walks under several types of lockdown. Our experiments show that completely closing medium and long distance travel to slow down the spread of a random walk is more efficient than than local restrictions.
\end{abstract}

\section{Introduction}\label{sec:intro} %%%%%%%%%%%%%%%%%%%%%%%%%%%
A common way to model real life phenomena relies on graphs. Here we describe a model where the edges of the graph represent possible routes a contaminant can take to travel, and vertices represents possible sites to infect. A simple random walk represents a theoretical immortal virus strain, and a self-replicating random walk would model an actual spread if one includes a replicating rate. In this work we aim to heuristically illustrate how the underlying geometry of travels plays a crucial role in the spread of a contaminant, as it is closely related to the \emph{conductance} of the graph (Definition \ref{conductance}), which governs the speed at which a simple random walk visits all the vertices. When the underlying geometry is flat, a random walk is slow, but in the presence of hyperbolicity, it becomes much faster.

\begin{figure}\label{cusps}
    \centering
    \includegraphics[width=\textwidth]{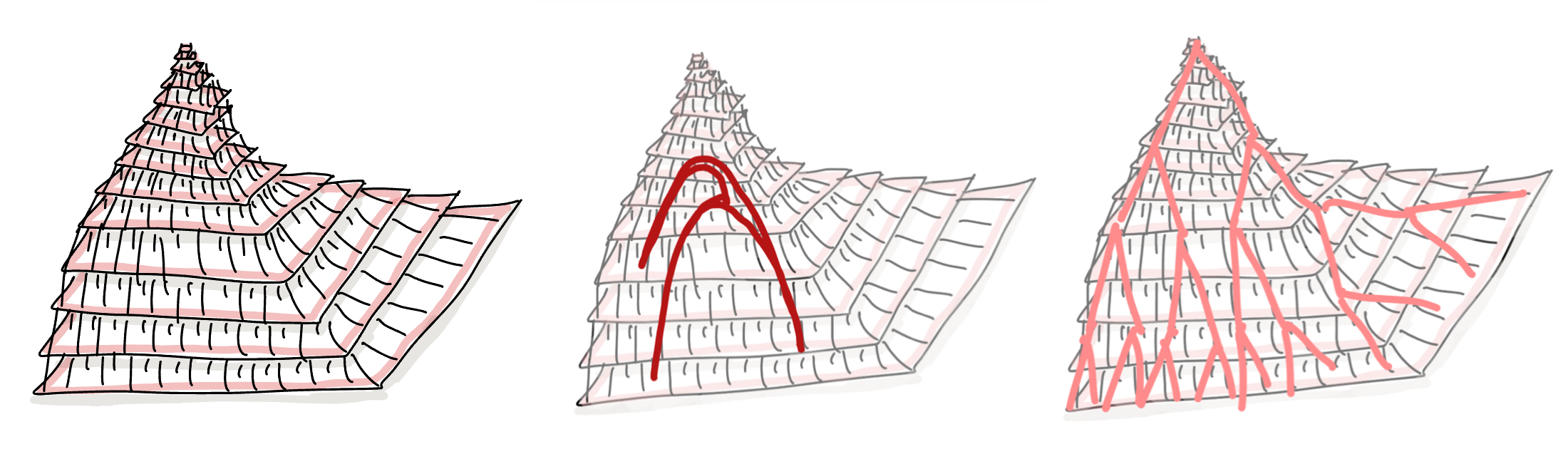}
    \caption{The first picture on the left shows the scaled down layers added. In the middle we added the behavior of geodesics traveling up and creating thin triangles, while on the right we see the underlying tree-like geometry.}
    \label{fig:RandomUniform20}
\end{figure}

Recall from the seminal work of Gromov in \cite{Gro} that a metric space (for instance a graph) is called {\it hyperbolic} if every triangle (that is, a triple of points and a choice of a shortest path, called {\it geodesic} between the three pairs of points) is uniformly thin regardless of the size of that triangle. Here, uniformly thin means that there is $\delta\geq 0$ such that for every triangle, every point on one of the 3 geodesics of that triangle is at distance less than $\delta$ from one of the two other geodesics. This is a tree-like geometry, that is most interesting when the diameter of the space is infinite, or large compared to $\delta$. According to \cite{Gro}\footnote{Example 1.8.A (b) on page 99} any metric space $(X,d)$ embeds as a horosphere in a hyperbolic metric space: a horosphere is a sphere centered at infinity. The construction consists in adding scaled down layers of the original space thus creating shortcuts in the metric space $(X,d)$, consisting in going up the layers to travel. With this construction, sets that had a small boundary in $X$, have a much larger one when embedded in the hyperbolic space, see Figure \ref{cusps}.
This construction allows to transform a flat geometry, with a slow random walk, into a hyperbolic one, that has much faster spreading random walks. In group theory, this is a well-known phenomena as groups with a flat geometry are in particular amenable, whereas non-elementary hyperbolic groups display a strong form of non-amenability.

We apply this idea to define a model of random graphs (Definition \ref{hrg}), that we call {\it horospherical random graphs}, defined by growing layers of graphs and whose limit is a particular case of {\it soft geometric random graphs}. Loosely speaking, we start with a metric space, choose vertices using some distribution, and then start at the first level by adding edges between points that are not too distant. Then at the second level, we add edges between points that are further away, etc. After a few levels, those graphs start to exhibit good expanding properties, in the sense that a random walk is very quick to visit all the vertices.

Modelling our random graphs according to the population density of a country (we worked on mainland France for our application in Section \ref{sec:application}), the edges represent possible routes that people can take, and hence the paths a virus can use to spread. Different layers of our sequence of graphs illustrate several confinement strategies, giving very different rates at which a random walk visits all the vertices. \\

\bigskip
%%%%%
{\bf Acknowledgements:} The authors thank Laurent Saloff-Coste for the reference \cite{LePeWi}, Goulnara Arzhantseva for discussions and references \cite{Pen}, \cite{Ham} and \cite{Nyb} as well as Patricia Reynaud-Bouret for interesting conversations and suggesting to look at the Small World features of horospherical graphs, and Damian Sawicki for pointing out a few inconsistencies in an earlier version. The authors also thank CNRS for making this collaboration materialize, and anonymous referees for interesting comments and suggestions that improved the paper.
%%%%
\section{Background on graphs}\label{sec:background} Recall that a {\it graph} $\Gamma$ is the data of a set $V$ of vertices and $E\subseteq V\times V$. Here we will assume our graphs to be undirected simplicial and connected.
The {\it adjacency matrix} $A=A_{(i,j)}$ of $\Gamma$ is an $m\times m$ matrix, where $m=|V|$ with values 1 when $\{i,j\}\in E$ and 0 otherwise. For an undirected simplicial graph, $A$ is a symmetric matrix with 0's on the diagonal. We further define the following
\begin{enumerate}
\item The {\it degree matrix} $D=D_{(i,j)}$ is a diagonal matrix with $D_{(i,i)}$ equals to the degree of the $i$-th vertex: $D_{(i,i)}=\sum_{j}A_{(i,j)}=d_i$ 
\item The {\it normalized adjacency matrix} is given by $\overline{A}=AD^{-1}$.
\item The {\it laplacian} is given by $L=D-A$.
\item The normalized laplacian is given by $\overline{L}=I-\sqrt{D^{-1}}A\sqrt{D^{-1}}$ (all its eigenvalues are in the interval $[0,1]$).
\end{enumerate}
\begin{defi}\label{conductance}For any two disjoint subsets of vertices $A,B\subseteq V$ , we denote $E(A,B)=E\cap(A\times B\cup B\times A)$, that is the set of edges of $\Gamma$ with one extremity in $A$ and the other in $B$. We denote by $\partial_E(A)=E(A,A^c)$, which are edges of the boundary of $A$. For a subset $A$ in $V$ we define its {\it volume} by $\nu(A)=\sum_{x\in A}\nu(x)$, where $\nu(x)=\left|\{y\in V \mid\{x,y\}\in E\}\right|$ is the {\it degree} of the vertex $x$. Hence $\nu(V)=\sum_{x\in V}\nu(x)$ is the total volume of the graph. The {\it conductance} of the graph $\Gamma$ is given by 
$$\varphi(\Gamma)=\min\left\{\frac{|\,\partial_E(S)\,|}{\min\{\nu(S),\nu(V\setminus S)\}}\ |\ \emptyset\not=S\subseteq V\right\}$$\end{defi}
The conductance is a measure on how well connected the graph $\Gamma$ is. This quantity is difficult to actually compute because it involves minimas over all subsets of a given set, but it is related to the first non-trivial eigenvalue of the normalized laplacian as follows (see for instance Theorem 14.4 of \cite{LePeWi} or Fact 6. page 5 of \cite{BuCh})
\begin{thm}[Cheeger inequality] Let $\Gamma$ be a connected graph and let $\lambda_1\in[0,1]$ be the smallest non-zero eigenvalue of its normalized laplacian. Then
$$\frac{\lambda_1}{2}\leq\varphi(\Gamma)\leq\sqrt{2\lambda_1}.$$
\end{thm}
\begin{rem}Many references (such as for instance \cite{Kow}) use the {\it expansion constant} (or sometimes called {\it Cheeger constant}), which is defined by 
$$h(\Gamma)=\min\left\{\frac{|\partial_E(A)|}{|A|}\hbox{ such that }\emptyset\not=A\subseteq V,\ 2|A|\leq|V|\right\}$$
with the convention that $h(\Gamma) = \infty$ if $\Gamma$ has at most one vertex. The expansion constant measures how difficult it is to separate the graph in two pieces of roughly the same size. The expansion constant is notoriously difficult to compute since the definition involves a minimum over all subsets. It is also related to the the smallest non-zero eigenvalue of the normalized laplacian, in a similar way as the conductance is.\end{rem}
We will be interested in simple random walks on connected graphs, a walk going to each neighbor with equal probability at each step.
\begin{defi}A {\it simple random walk} on $\Gamma$ is a sequence $(X_n)_{n\geq 0}$ of $V$-valued random variables (measurable functions) defined on a common probability space $(\Omega,\Sigma,P)$, such that, for any $n\geq 0$, and any vertices $x_0,\dots, x_n$ and $y\in V$, whenever $\{x_i,x_{i+1}\}\in E$, the following holds
$$P\{X_{n+1} = y \ |\  (X_n,\dots,X_0) = (x_n,\dots,x_0)\} = P\{X_{n+1} = y \ |\ X_n = x_n\}$$
with the conditional probability of the walk being at $y\in V$ at time $n+1$ given that it is at $x_n$ at time $n$ given by
$$P\{X_{n+1} = y \ |\  X_n = x_n\}=\left\{\begin{array}{cc}0 &\hbox{ if }\{x_n,y\}\not\in E\\
\frac{1}{|\{x\in V |\{x_n,x\}\in E\}|} & \hbox{ if }\{x_n,y\}\in E\end{array}\right.$$
This means that the walk at time $n+1$ doesn't depend on what happened on the $n-1$ first steps, only on its position in the graph at time $n$.
\end{defi}
If the graph is not bipartite, then the walker is, after a time $n$ large enough, almost as likely to be located at any of the vertices of the graph, independently of the starting point or of the steps taken to get there. For a particle as likely to be anywhere on the graph, the probability governing its position is the uniform probability, but for a particle starting somewhere on the graph, this probability is never exactly uniform, but almost. This means that {\it the walk converges to the uniform distribution}, and the speed at which this convergence occurs, namely the number $n$ needed for the walker's distribution to be close to uniform is related to the first eigenvalue $\lambda_1$ as follows (see \cite{BuCh} or \cite{LePeWi})
$$\left|P(X_n=x)-\frac{\nu(x)}{N}\right|\leq\sqrt{\frac{\nu(x)}{v_-}}|\lambda_1-1|^n$$
where $N=\nu(V)$ is the total volume of the graph and $v_-$ is the lowest degree of a vertex. Intuitively, this means that the walk could be anywhere in the graph, hence has been everywhere in the graph. The value $\lambda_1$ is often called {\it spectral gap}, and the larger it is, the faster the sequence $|\lambda_1-1|^n$ tends to 0 (since $\lambda_1\in[0,1])$.
\begin{defi}The {\it clustering coefficient} at a vertex $v_i\in V$ is given by 
$$C_i=2\sum_{j,k}\frac{A_{(i,j)}A_{(j,k)}A_{(k,i)}}{d_i(d_i-1)}$$
and the {\it average clustering coefficient} is the mean value over all vertices, of local clustering coefficients. 
\end{defi}
\begin{rem} The clustering coefficient at a vertex counts the proportion of closed paths of length 3, among all neighbors of that vertex. Clustering coefficients are used to give a measure of {\it Small World}, which are sparse graphs with a large number of nodes, a large average clustering coefficient and small average distances between pairs of points in the graph, see \cite{Dur} or \cite{Wat}.\end{rem}
\section{Description of the probabilistic graph}\label{sec:graph_model} 
The general settings for constructing our horospherical random graphs is a particular case of soft random geometric graphs described by Penrose \cite{Pen}, and is described in the following. For a metric space $(X,d)$, $x\in X$ and $r\geq 0$ we denote by 
$$B_x(r)=\{y\in X\ |\ d(x,y)\leq r\}$$
the {\it ball of radius $r$ centered at $x$}.
\begin{defi}[Horospherical random graphs]\label{hrg}
Suppose that $(X, d)$ is a discrete metric space, locally finite (balls of finite radius have finitely many elements, for instance obtained from a finite subset of $\R^n$ with the induced Euclidean metric from $\R^n$). Moreover, suppose that $0<r_1<r_2<\ldots$ is an increasing sequence of positive real numbers and $\phi=(p_1,p_2\ldots)$, $p_i\in[0,1]$ is any sequence of real numbers in the unit interval. This defines a piece-wise constant function $\phi$, taking the value $p_k$ in the interval $[r_k,r_{k+1})$, for any $k\in\N$, called {\it connection function}. We define a collection of random graphs $\mathcal{G} = \{\Gamma_n\}_{n=1}^\infty$ in the following way
\begin{itemize}
%    \item $\Gamma_0 = (V_0,E_0)$ is the standard square grid with vertices at pairs of integer points.
    \item $\Gamma_1 = (X,E_1)$ where for all $x\in X$ and for all $x'\in B_x(r_1)\subseteq V$, the edge $\{x,x'\}$ is added to $E_1$ with probability $p_1$.
    \item $\Gamma_n = (X,E_n)$ where the edge set is generated randomly as follows:
    \begin{itemize}
        \item $E_{n-1}\subseteq E_n$
        \item For all $x\in X$ and for all $x'\in X\cap \left(B_x(r_n)\setminus B_x(r_{n-1})\right)$, the edge $\{x,x'\}$ is added to $E_n$ with probability $p_n$.
    \end{itemize}
\end{itemize}
The collection $\mathcal{G}$ is what we call {\it horospherical random graph}, and for each $i\in{\bf N}$, the graph $\Gamma_i$ is the {\it $i$-th layer}, or the {\it graph at height $i$}. If the connection function $\phi$ is finitely supported, then the collection stabilizes at $\Gamma_n$ for $n$ large enough.
\end{defi}
%%%
\begin{rem}
If $\Gamma_1$ is connected, then $\Gamma_n$ is connected for all $n>1$. Moreover, from the definition of Cheeger constant it follows that $h(\Gamma_i)\leq h(\Gamma_j)$ when $i\leq j$ and according to \cite{CvDoSa} Theorem 2.1 $\lambda_k(\Gamma_j)\geq\lambda_k(\Gamma_i)$. Our goal will be to study the behavior of $\lambda_1(\Gamma_i)$ as $i$ increases. 
\end{rem}
%%%
\begin{rem}Denote by
$$W_r = \{(x,y)\in X\times X\mid d(x,y) = r\}$$
the cylinder of diameter $r$ around the diagonal in $X\times X$. Our sequence of horospherical random graphs is a particular case of the more general following setting, where we start with $(X,d)$ a metric space, with the following additional data.
\begin{itemize} 
    \item Let $\phi:\R^+\to [0,1]$ be a measurable function, called {\it connection function} (the above Definition \ref{hrg} is the case where this function is piecewise constant).
    \item  Define $\hat{W}_r\subseteq W_r$ to be a random subset where each $w\in\hat{W}_r$ appears with probability $\phi(r)$ (this is ($|W_r|$) samples from a Bernoulli($\phi(r)$) distribution).
    \item Define $E_r = \bigcup_{r'\le r}\hat{W}_r$ (we can also view $\hat{W}:\R^+\to X\times X$
        where $r\mapsto \hat{W}_r$).
    \item The graph $\Gamma_r = (X,E_r)$ is the graph obtained at height $r$.
\end{itemize}
For modelling purpose and to recover our definition and several models from the literature, we restrict our attention to a discretization of the above construction, done as follows. Choose $h>0$, which will be the step size and define $\bar{\phi}:\R^+\to[0,1]$ to be the step (i.e. piece-wise constant) function obtained from $p$ above on the intervals $[kh,(k+1)h)$, where $k\in\N$ by 
$$\frac{1}{h}\int_{kh}^{(k+1)h}\phi(t)dt$$ 
and we recover the $i$-th layer $\Gamma_i=E_{r_i}$ from Definition \ref{hrg}. 
\end{rem}
\begin{rem}
    In Definition \ref{hrg} the sequence $(r_i)$ doesn't need to go to infinity, but if it is bounded, uniform local finiteness will force the differences $B_x(r_n)\setminus B_x(r_{n-1})$ to eventually be empty and the sequence of graphs stabilizes. Similarly, when the metric space $X$ has bounded diameter, the sequence of graphs will stabilize regardless of the sequence $(r_i)$ going to infinity.
\end{rem}
\section{Related models}
Let us first notice that for a finitely supported connection function $\phi$, or for a finite diameter metric space $(X,d)$ the sequence of graphs stabilizes after some height.

The random Erd\H{o}s-R\'enyi graphs are instances of horospherical random graphs in case where the parameter $p$ is a constant function. The first level of a horospherical random graph is a geometric random graph, see \cite{Ham} or \cite{Nyb} for related studies. The limiting graph in our horospherical random graphs are instances of soft geometric random graphs studied by Penrose \cite{Pen}. For us, keeping track of the layers is a way of understanding how the long distance connections influence the behavior of the first eigenvalue and hence of the conductance.

When $X\subset\R^d$ our horospherical random graphs $\Gamma_n$ are instances of geometric inhomogeneous random graphs studied by \cite{BrKeLe}, also called spatial inhomogeneous random graphs in \cite{HoHoMa}, in the particular case where all the weights at vertices are constants. If the metric space one starts with, is the hyperbolic plane in the Poincaré disk model, horospherical random graphs are generalizations of hyperbolic random graphs: we do not require the underlying metric space to be hyperbolic, but choosing the connection function well, the resulting graph could be hyperbolic. The connection function for hyperbolic random graphs is studied in \cite{GuPaPe}.
\section{Application}\label{sec:application}
By using the first graphs in the sequence of graphs defined in Section~\ref{sec:graph_model}, we can produce a model that simulates the effects of lockdown measures taken during a pandemic. Our model has two sets of parameters, a nondecreasing sequence of distances $0\le r_1\le r_2\le \ldots$ and a sequence of probabilities $p_1,p_2,\ldots\in[0,1]$, giving a connection function $\phi$. The distance $r_i$ represents the possible range of interaction and the probability $p_i$ represents the probability of interaction between the vertices. In this section, we produce simulations for two different sets of vertices and examine the resulting $\lambda_1$ values\footnote{Code for this project available at: \url{https://gist.github.com/azlawson/a545b341f5621cabbb649e75a55a4fca}}.
\subsection{Uniformly distributed vertices}
For our first simulation, we will explore some statistical quantities for horospherical random graphs with five different probability parameter functions $p_i$ for $i=1,\ldots, 5$. In this case, we generate a single random set of 200 points sampled uniformly from an 8 by 8 square (those numbers were chosen after trial and error to get the clearer pictures in the simulations of Subsection \ref{France}). For each probability parameter and each point set, we generate 20 horospherical random graphs by adding a level at each of 100 evenly spaced values of $r$ between 0 and $8\sqrt{2}$. For each level we recorded the values of $\lambda_1$, maximum degree, average degree, and sparsity (that is, the proportion of edges compared to the maximal possible number of edges, which is $m(m-1)/2$ if $m$ is the total number of vertices). Our goal is to explore the similarities in the different choices of $p_i$ within and across the different point sets as well as the variance in these selections. Figure \ref{fig:RandomUniform20} contains four plots, one for each of these statistics. Each plot contains the average value considered along with a band of one standard deviation for each probability parameter function. These functions consist of two constant functions, $\phi_1(r) = 0.5$ and $\phi_2(r)=1$, an affine function $\phi_3(r) = \frac{-r}{8\sqrt{2}}+1$, an exponential $\phi_4(r) = e^{-r}$, and $\phi_5(r) = \frac{\exp\{-(r-\sqrt{32})^2/16\}}{\sqrt{16\pi}}$, a normal density with mean $\sqrt{32}$ and standard deviation $2$.
\begin{figure}
    \centering
    \includegraphics[width=\textwidth]{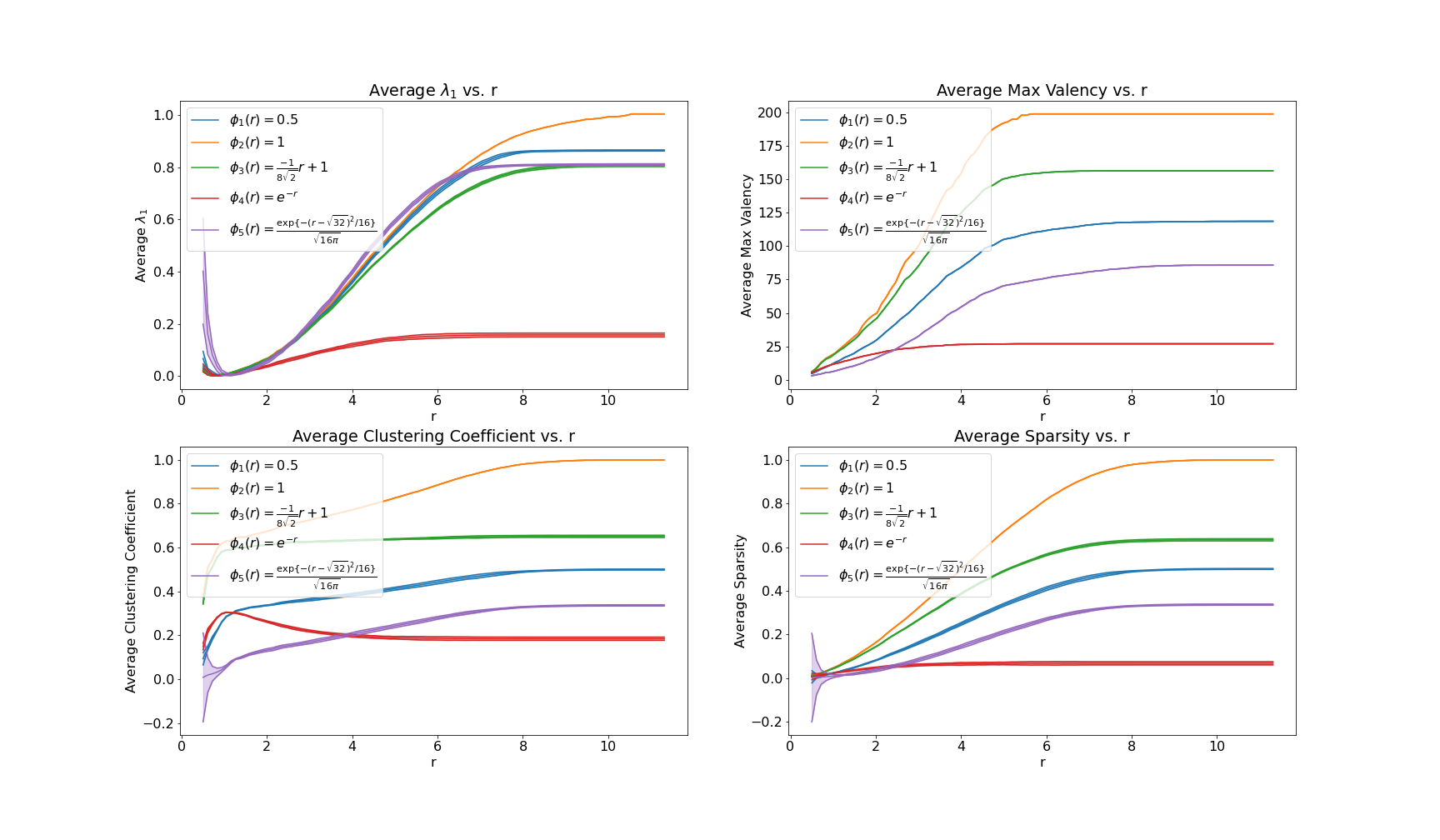}
    \caption{Four plots of graph properties over 5 different connection functions. For each connection function and each property, the given property was averaged over 20 iterations on a uniformly distributed set, then plotted with a band of one standard deviation.}
    \label{fig:RandomUniform20}
\end{figure}
\subsection{A model for mainland France}\label{France}
\begin{figure}
    \centering
    \includegraphics[width=\textwidth]{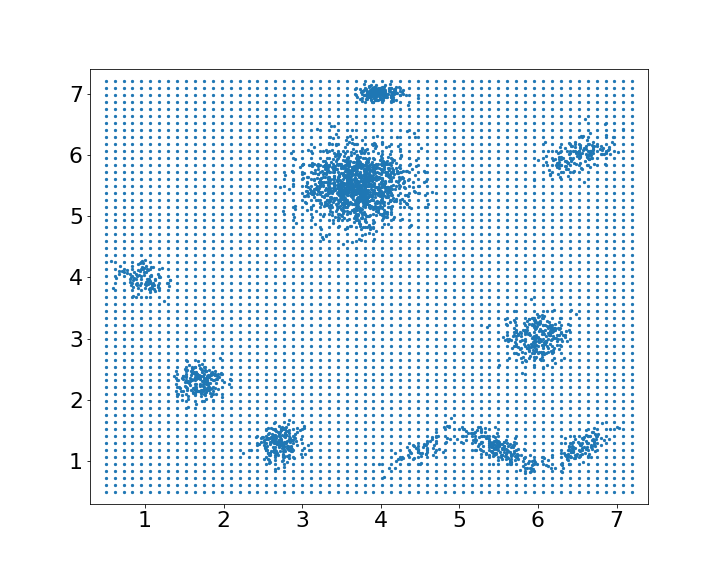}
    \caption{The vertex set $X$, obtained with a simulation of major metropolitan areas of mainland France along with a grid of points representing the rest of the population.}
    \label{fig:Franceverts}
\end{figure}

We construct a model of mainland France by selecting the most populated metropolitan areas such as Paris, Toulouse, Bordeaux, Lyon, Strasbourg, Nantes, Lille, and C\^ote d'Azur (comprised of Nice, Marseille, and Toulon), using the picture from Wikipedia\footnote{https://fr.wikipedia.org/wiki/Liste$\underline{\ }$des$\underline{\ }$aires$\underline{\ }$urbaines$\underline{\ }$de$\underline{\ }$France} copied on an $8\times 8$ square. The vertices of our graph sequence were drawn from bivariate normal distributions centered at each of those cities. The number of points drawn from each city is proportional to the population of the city relative to the population of mainland France. Finally, we overlayed a grid of 60 by 60 points (so a mesh of $8/60<0.15$) to represent the distribution of the remaining population throughout the country side, and will allow us to force all graphs to be connected. Figure \ref{fig:Franceverts} depicts the resulting vertices, 6'000 dots, so that each point is then roughly 10'000 people if we estimate mainland France population by 60 million.
On this set of vertices, denoted by $X$, we compute several simulations with four different connection functions, named $U, C, S$ and $I$, heuristically depicting daily people movements in different situations, and that we now describe. First, we consider the connection function given by
\[U(r) = \begin{cases}
1 & r\le 0.15\\
0.05 & 0.15<r\le0.3\\
0.03 & 0.3<r\le 1\\
0.02 & 1<r\le3\\
0.01 & 3<r\le 10\\
\end{cases}
\]
with the corresponding random graph $\Gamma_U = (X,E_{10})$. This corresponds to a daily situation that we estimated very roughly using 2019 official data available on-line \footnote{https://www.ecologie.gouv.fr/statistiques-du-trafic-aerien and\\ https://www.autorite-transports.fr/wp-content/uploads/2021/01/bilan-ferroviaire-2019.pdf} before any type of lockdown. Next, we consider the connection function given by
\[S(r) = \begin{cases}
1 & r\le 0.15\\
0.1 & 0.15<r\le0.3\\
0.05 & 0.3<r\le 1\\
0.01 & 1<r\le3\\
0.005 & 3<r\le 10\\
\end{cases}
\]
with corresponding graph $\Gamma_S = (X,E_{10})$. This corresponds to some soft travel restrictions, simulating the case where all means of travel are open, but are operating at limited capacity. Here there is an increase of short travels compared to the pre-lockdown situation, that in fact has very low impact on the first eigenvalue. The connection function
\[C(r) = \begin{cases}
1 & r\le 0.15\\
0.05 & 0.15<r\le0.3\\
0.003 & 0.3<r\le 1\\
0.002 & 1<r\le3\\
0.001 & 3<r\le 10\\
\end{cases}
\]
corresponds to our estimates of long distance travels in France during the 2020 first lockdown, one of the hardest ones, where some 90 percent of airplane activity was cancelled, with the corresponding random graph $\Gamma_C = (X,E_{10})$. Finally the connection function
\[I(r) = \begin{cases}1 & r\le 0.3\\0&\hbox{otherwise}
\end{cases}
\]
with corresponding graph $\Gamma_I = (X,E_{10})$. This corresponds to a theoretical lockdown, where long distance travel is completely stopped but locally, everything is widely open, without restrictions. In that way, travel between cities is eventually possible, but takes place at a very slow pace.

\medskip

We seek to compare these graphs through various statistics. Specifically, for each graph we compute the first eigenvalue $\lambda_1$, Sparsity, Max degree, and Average degree.  Because $X$ is fixed and $p_1\equiv 1$, the graph $\Gamma_I$ is not random, and the probability 1 for $r_1\leq 0.15$ on all the graphs is to ensure that those graphs are connected. However, the next layers in the graphs $\Gamma_U, \Gamma_S$ and $\Gamma_C$ are random; thus, we choose to compute 20 iterations and average the relevant statistics. The results are summarized in Table~\ref{tab:Exp1Stats}. We can see that the graphs are nearly equal in Sparsity, except for $\Gamma_C$, the hard confinement graph which is much more sparse than the three others and has significantly lower max and average  degree while its first eigenvalue $\lambda_1$ is roughly a third of the ones in the unconfined graph $\Gamma_U$ and the lightly confined model $\Gamma_S$. But comparing it to $\Gamma_I$, we notice that $\lambda_1(\Gamma_I)$ is $3.5\times 10^{-4}$, so several orders of magnitude smaller than the other graphs. This indicates that random walks on $\Gamma_I$ (and unsurprisingly on, $\Gamma_U$ and $\Gamma_S$)  will spread much faster than random walks on $\Gamma_I$.  

\medskip

\begin{figure}
    \centering
    \includegraphics[width=\textwidth]{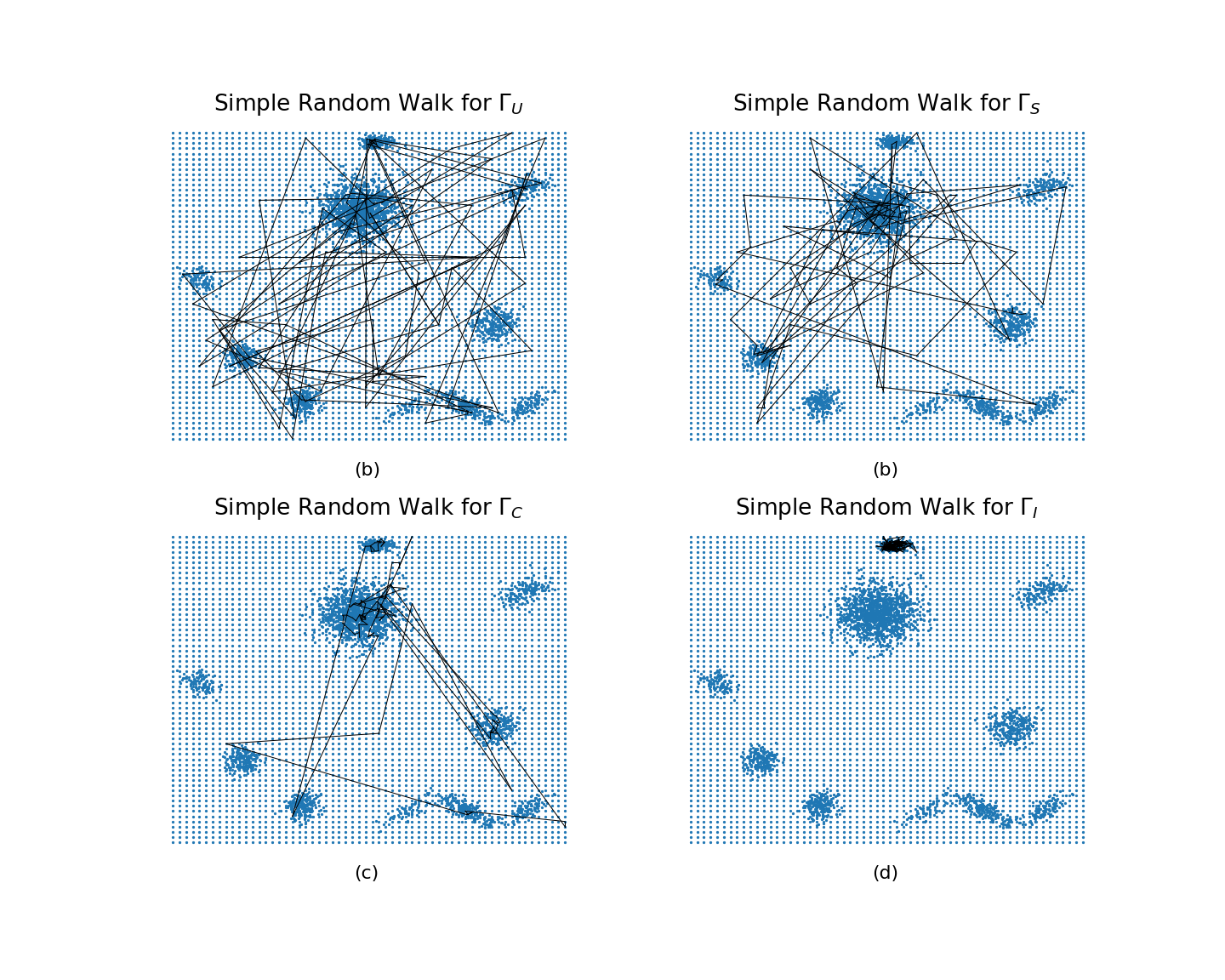}
    \caption{A simple random walk on $\Gamma_U, \Gamma_S, \Gamma_C$ and $\Gamma_I$.}
    \label{fig:Exp1SimpleWalk}
\end{figure}

To illustrate this, we consider two different types of random walks. First, we simulate a simple random walk on the four graphs. For each graph, we generate 100 different 100-step random walks with a random starting points. For each walk, we compute the average and maximum pairwise distances. These numbers suggest that traveling long distances is much easier to achieve in the set up for $\Gamma_C$ than for that of $\Gamma_I$. We summarize these numbers in Table~\ref{tab:Exp1WalkStats}.
Figure~\ref{fig:Exp1SimpleWalk} illustrates an example walk resulting from the simulation. We can see that, as expected, the walk for $\Gamma_U$ very quickly spreads everywhere, while the walk for $\Gamma_I$ remains concentrated near its inception. What is interesting is the fact that for $\Gamma_S$ and $\Gamma_C$, the walk spreads pretty much everywhere as well.
\begin{table}
    \centering
    \bgroup
    \def\arraystretch{1.25}%
    \begin{tabular}{|c|c|c|c|c|}\hline
            & $\Gamma_U$& $\Gamma_S$ &$\Gamma_C$ & $\Gamma_I$  \\\hline
$\lambda_1$ & 0.46 & 0.29 & 0.09 & 0.00035\\\hline 
        Sparsity & 0.019 & 0.016 & 0.006 & 0.015\\\hline
        Max  degree & 280&278&178&470 \\\hline
        Average  degree &121&96 &40 &97\\\hline
        Clustering Coefficient & 0.06 & 0.08 & 0.19 & 0.62\\\hline
    \end{tabular}
    \egroup
    \caption{Various statistics for $\Gamma_U,\Gamma_S,\Gamma_C$ and $\Gamma_I$.}
    \label{tab:Exp1Stats}
\end{table}
\begin{table}
    \centering
    \bgroup
    \def\arraystretch{1.25}%
    \begin{tabular}{|c|c|c|c|c|}\hline
        & $\Gamma_U$ & $\Gamma_S$ & $\Gamma_C$ & $\Gamma_I$  \\\hline
        Mean of mean distances & 3.07&2.87&2.57&0.27\\\hline 
        St. dev. of mean distances &0.23&0.35&0.63&0.20\\\hline
        Mean of max distances & 8.2& 7.84&7.2&0.91\\\hline
        St. dev. of max distances&0.43&0.53&0.64&0.49\\\hline
    \end{tabular}
    \egroup
    \caption{Various statistics for 100 walks on $\Gamma_U,\Gamma_S,\Gamma_C$ and $\Gamma_I$}
    \label{tab:Exp1WalkStats}
\end{table}

\begin{figure}
    \centering
    \includegraphics[width=\textwidth]{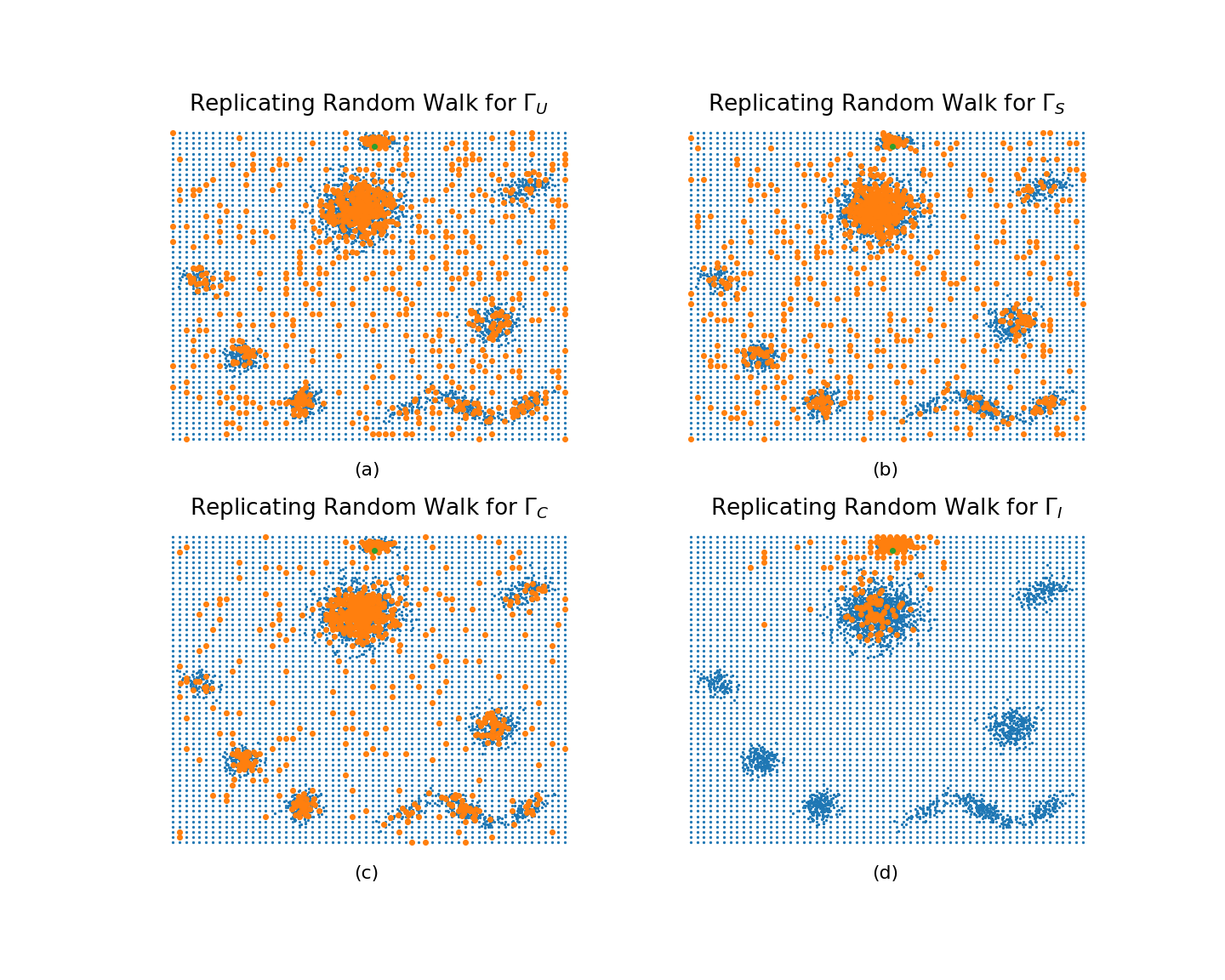}
    \caption{A replicating random walk simulation on $\Gamma_U, \Gamma_S, \Gamma_C$ and $\Gamma_I$.}
    \label{fig:Exp1RepWalk}
\end{figure}
For another illustration, we compute a replicating random walk. The walk for each graph begins with one particle, placed on a node $\rho_1$ in the same starting location (near Lille, on the top of the grid). After 10 steps of the random walk, if the particle is on a node $\rho_2\not=\rho_1$, then we duplicate the particle. We now have two particles that embark on two simultaneous random walks. After 10 more steps, any particle on a node that has not been the place of a duplication before, will duplicate. The process continues and this walk proceeds until the first particle has walked 100 steps. Figure~\ref{fig:Exp1RepWalk} shows the initial positions of each particle at the end of the replicating walk. Again, the difference in spread between the different graphs are clear.

\medskip

\begin{figure}
    \centering
    \includegraphics[width=0.5\textwidth]{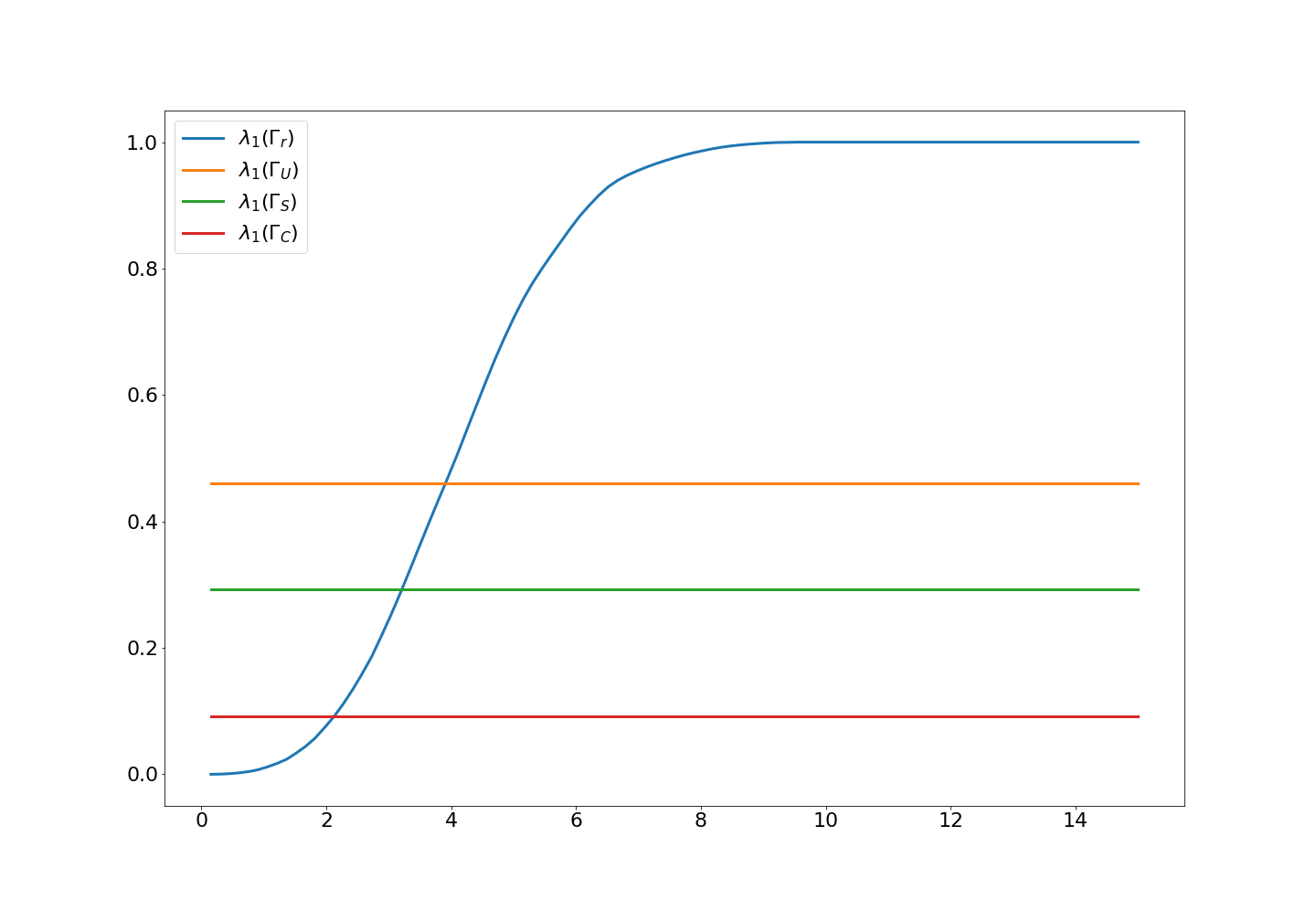}
    \caption{A graph of $\lambda_1(\Gamma(r))$. The horizontal lines are $\lambda_1(\Gamma_U), \lambda_1(\Gamma_S)$ and $\lambda_1(\Gamma_C)$.}
    \label{fig:Prob1Dist}
\end{figure}

Finally, we explore the effect of the distance parameter, $r$ when $p\equiv 1$, which we will denote as $\Gamma(r)$. Clearly, as $r$ increases, $\lambda_1(\Gamma(r))$ will increase as well. Figure~\ref{fig:Prob1Dist} shows the behavior of $\lambda_1(\Gamma(r))$ in relation to $r$. The orange line in the Figure is $\lambda_1(\Gamma_U)$, the green line is $\lambda_1(\Gamma_S)$ and the red line is $\lambda_1(\Gamma_C)$, the lowest eigenvalue. It is interesting to note that the $\lambda_1(\Gamma(r)) = \lambda_1(\Gamma_C)$ when $r\approx 2.3$. In our model, this equates to a distance of about 300 km, or 180 miles.
\section{Discussion on the findings and future work}
The presence of edges between far away points dramatically increases the speed at which the random walk visits all the vertices of the graph. From the view point of the random walks, completely shutting down medium and long distance travel (which amounts to looking at low level graphs) is more efficient in slowing down the spread of the walk than restricting local interactions. In fact even with high local clustering coefficient, corresponding to a level one graph with large edge probability, the random walk is much slower than in a sparse graph with edges linking far away points. 

Our heuristic findings show that completely cutting travels over 300 kilometers would give a random walk much slower to spread than in a graph with strict local confinement, amounting to low average  degree. Of course, at small scale a high local clustering coefficient would give a high local spread once the walk has reached a location, and this is what strict lockdowns try to mitigate.

Low level graphs, with small first eigenvalue $\lambda_1$, also have small conductance and are easy to cut into pieces of similar sizes, meaning that clusters could be disconnected from the rest of the graph. In terms of random walks, this partial confinement mitigates the damages: if the random walk is slow enough, one can divide a graph in $n$ pieces $R_1,\dots,R_n$ easily disconnected by cutting relatively few edges. Compared to their sizes, the $R_i$'s have {\it small isoperimetric inequality}, meaning that to isolate a piece $R_i$ of the rest of the graph to trap the random walk there can be done removing relatively few edges. This can be seen in the graph $\Gamma_I$, where the replicating walk starting in Lille doesn't spread beyond Paris, whereas on $\Gamma_U,\Gamma_S$ and even on $\Gamma_C$ it has spread everywhere.

The impact of long distance travel in the spread of diseases is probably well-known by specialists, and the mathematics described in this paper are very classical, but the conclusions seem to be controversial\footnote{https://www.icao.int/Newsroom/Pages/Update-on-ICAO-and-WHO-Coronavirus-Recommendations.aspx}, even if a scientific consensus along those lines seems to emerge \cite{Adiga2020}.

Obviously, comparing our model with actual data, making finer estimates and more accurate models would be interesting and could provide interesting information, but beyond our computing capacity and this article aims to explain the mathematical phenomena that makes long distance travel - even very little - drastically speed the spread of a random walk.
\bibliographystyle{plain}
\bibliography{main.bib}
\end{document}